\documentstyle[11pt]{article}
\newlength{\standardunitlength}
 \newtheorem{lemma}{Lemma}
\newtheorem{theorem}{Theorem} 
\newenvironment{proof}{\noindent {\sc Proof:}}{$\Box$ \vspace{2 ex}}

\begin{document}

\begin{center}
Semisimple orbits of Lie algebras and card-shuffling measures on Coxeter groups
\end{center}

\begin{center}
By Jason Fulman
\end{center}

\begin{center}
Dartmouth College
\end{center}

\begin{center}
Department of Mathematics
\end{center}

\begin{center}
6188 Bradley Hall
\end{center}

\begin{center}
Hanover, NH 03755, USA
\end{center}

\begin{center}
email:jason.e.fulman@dartmouth.edu
\end{center}

\begin{center}
\end{center}

\begin{center}
\end{center}

\begin{center}
\end{center}

\begin{center}
\end{center}

\begin{center}
\end{center}

\begin{center}
\end{center}

\begin{center}
\end{center}

\begin{center}
\end{center}

1991 AMS Primary Subject Classifications: 20G40, 20F55

\newpage
\begin{abstract}
	Solomon's descent algebra is used to define a family of signed measures $M_{W,x}$
for a finite Coxeter group $W$ and $x \neq 0$. It is known that the measures
corresponding to $W$ of types $A$ and $B$ arise from the theory of card shuffling
and are related to the Poincare-Birkhoff-Witt theorem and splitting of Hochschild
homology. Formulas for these measures and obtained in special cases. The
eigenvalues of the associated Markov chains are computed. By elementary algebraic
group theory, choosing a random semisimple orbit on a Lie algebra corresponding
to a finite group of Lie type $G^F$ induces a measure on the conjugacy classes of
the Weyl group $W$ of
$G^F$. It is conjectured that this measure on conjugacy classes is equal to the
measure arising from $M_{W,q}$ (and further that
$M_{W,q}$ is non-negative on all elements of $W$). This conjecture is proved
for all types for the identity conjugacy class of $W$, and is confirmed for all
conjugacy classes for types
$A_n$ and $B_n$.
\end{abstract}    

\section{Definition of the Signed Measures $M_{W,x}$} \label{definition}

	This section defines signed measures $M_{W,x}$ for any Coxeter group $W$ and
real $x \neq 0$. By a signed measure is meant an element of the group
algebra $Q[W]$ of $W$ whose coefficients sum to one. The motivation	for this
definition comes from work of Bergeron, Bergeron, Howlett, Taylor \cite{BBHT} and
Bergeron and Bergeron \cite{BB}. For types $A$ and $B$, results of Bergeron and
Wolfgang \cite{BergeronWolfgang} show that $M_{W,x}$ is related to the
Poincare-Birkhoff-Witt theorem and splitting of Hoschild homology.

	Let $\Pi$ be a set of fundamental roots for a root system of $W$. Call subsets
$K_1$ and $K_2$ of $\Pi$ equivalent if there is a $w$ such that $w(K_1)=K_2$. Let
$\lambda$ be an equivalence class of subsets of $\Pi$ under the action of $W$ and
let $\lambda_K$ be the equivalence class of the set $K$. Let $|\lambda|$ denote
the size of the equivalence class $\lambda$, and let $\| \lambda \|$ denote the
size of the set
$K$ for any $K \in \lambda$. 

 For $w \in W$, define $D(w)$ as the set of simple positive roots mapped to
negative roots by $w$ (also called the descent set of $w$). Let $d(w)=|D(w)|$.
For $J \subseteq \Pi$, let
$X_J=\{w
\in W| D(w) \cap J = \emptyset \}$ and $x_J = \sum_{w \in X_J} w$. For $K
\subseteq J
\subseteq \Pi$ define $\mu_K^J = \frac{|\{w \in X_J: w(K) \subseteq \Pi
\}|}{|\lambda_K|}$. Set $\mu_K^J=0$ if $K \not \subseteq J$. Since the matrix
$(\mu_K^J)$ is upper triangular with non-zero diagonal entries, it is invertible.
Letting
$(\beta_K^J)$ be its inverse, define
$e_J$ and $e_{\lambda}$ in the descent algebra of $W$ by

\begin{eqnarray*}
e_J & = & \sum_{K \subseteq J} \beta_K^J x_K\\
e_{\lambda} & = & \sum_{J \in \lambda} \frac{e_J}{|\lambda|}
\end{eqnarray*}

	Bergeron, Bergeron, Howlett, and Taylor \cite{BBHT} prove that the $e_{\lambda}$
are orthogonal idempotents of the descent algebra decomposing the identity.

{\bf Definition} For $W$ a finite Coxeter group and $x \neq 0$,
define a signed measure $M_{W,x}$ on $W$ by

\[ M_{W,x} = \sum_{\lambda} \frac{e_{\lambda}}{x^{\| \lambda \|}} \]

	For $w \in W$, let $M_{W,x}(w)$ be the coefficient of $w$ in $M_{W,x}$.

\begin{theorem} \label{signedmeasure} $M_{W,x}$ is a signed measure on $W$.
\end{theorem}

\begin{proof}
	Writing each $e_{\lambda}$ as $\sum_{w \in W} c_{\lambda}(w) w$ it must be
proved that 

\[ \sum_{w,\lambda} \frac{c_{\lambda}(w)}{x^{\| \lambda \|}} = 1 \]

	This clearly follows from the stronger assertion that:

\[ \sum_{w} c_{\lambda}(w) = \left\{ \begin{array}{ll}
																																				0 & \mbox{if $\| \lambda \|>0$}\\
																																				1 & \mbox{if $\| \lambda \|=0$}
																																				\end{array}
						\right. \]

	Corollary 6.7 of Bergeron, Bergeron, Howlett, and Taylor \cite{BBHT} implies that
$e_{\emptyset} =
\frac{\sum_{w
\in W} w}{|W|}$. Thus $\sum_w c_{\lambda}(w) = 1$ if $\|
\lambda
\|=0$. Since the
$e_{\lambda}$ are idempotents, the value of $\sum_{w} c_{\lambda}(w)$ is
either 0 or 1. Since $\sum_{\lambda} e_{\lambda}=1$, clearly $\sum_{w,\lambda}
c_{\lambda}(w)=1$. Combining this with the fact that $\sum_w c_{\lambda}(w) = 1$
if $\| \lambda \|=0$ shows that $\sum_w c_{\lambda}(w) = 0$ if $\| \lambda \|>0$.
\end{proof}	

\vspace{.5mm}

{\bf Remarks}

\begin{enumerate}

\item If $W=S_n$, then as noted in Bergeron and Bergeron \cite{BB}, the 
measure $M_{W,x}$ corresponds to performing an
$x$-shuffle on $W$ according to the Gilbert-Shannon-Reeds model of card
shuffling. This model of card shuffling is described clearly and analyzed by
Bayer and Diaconis \cite{BD}. Let $d(w)=|D(w)|$. Bayer and Diaconis prove
combinatorially that

\[ M_{S_n,x}(w) = \frac{{x+n-1-d(w) \choose n}}{x^n} \]

	Some further information about the measure $M_{S_n,x}$ can be found in
Fulman \cite{F}. For instance a generating function is derived for
the distribution of the length of a permutation (in terms of the generators
$\{(1,2),(2,3),\cdots,(n-1,n)\}$) chosen from this measure.	

	For $W$ of type $B$ (and thus also of type $C$), Bergeron and Bergeron
\cite{BB} prove that

\[ M_{B_n,x}(w) = \frac{(x+2n-1-2d(w))(x+2n-3-2d(w)) \cdots (x+1-2d(w))}{x^n n!}
\]

	An easy computation using formulas at the end of Section 2 of Bergeron and
Bergeron \cite{BB} proves that

\[ M_{I_2(p),x}(w) =             \left\{ \begin{array}{ll}
																																				\frac{(x+1)(x+p-1)}{2px^2} & \mbox{if
$d(w)=0$}\\
																																				\frac{(x+1)(x-1)}{2px^2} & \mbox{if
$d(w)=1$}\\
																																				\frac{(x-1)(x-p+1)}{2px^2}		& \mbox{if
$d(w)=2$}
																																				\end{array}
			\right.			 \]

	From the definition of $M_{W,x}$, it is clear that $M_{W,x}(w)$ depends only
on $D(w)$, the descent set of $w$. Results and conjectures for other $W$ appear
in Section \ref{eigenvalues}.

\item Observe that the $x \rightarrow \infty$ limit of $M_{W,x}$ is the uniform
distribution on $W$ (this follows from the formula for $e_{\emptyset}$ in the
proof of Theorem \ref{signedmeasure}). The eigenvalue computations of Section
\ref{eigenvalues} can be used to give results on how fast this convergence occurs.

\item The elements $M_{W,x}$ of the group algebra of $W$ convolve nicely in the
sense that $M_{W,x} M_{W,y} = M_{W,xy}$. This follows from the fact that the
$e_{\lambda}$ are orthogonal idempotents.

\item As will emerge, $M_{W,x}(w)$ need not always be positive. Part of Conjecture
1 of Section \ref{semisimple} states that $M_{W,q}(w) \geq 0$ if
$W$ is a Weyl group of a finite group of Lie type and $q$ is a power of a prime
which is regular and good for $W$ (these terms are defined in Section
\ref{semisimple}).

\item Bidigare, Hanlon, and Rockmore \cite{BHR} define and study interesting
random walks on the chambers of hyperplane arrangements. Bidigare
\cite{B} defines a face algebra associated to a hyperplane arrangement
and shows that if the hyperplane arrangement comes from a reflection group $W$,
then the descent algebra of $W$ is anti-isomorphic to the trivial isotypic
subalgebra of the face algebra. This suggests that the measures $M_{W,x}$ are
special cases of the Bidigare-Hanlon-Rockmore measures. This is known to be true
for $W$ of type $A$.

\end{enumerate}

\section{Formulas for $M_{W,x}$ and the Eigenvalues of the Markov Chain
Associated to $M_{W,x}$} \label{eigenvalues}

	This section considers formulas for $M_{W,x}$. A expression is found for
$M_{G_2,x}$, and for all $W$, the values of $M_{W,x}$ on the identity and longest
element of $W$ are computed. This will allow us compute the eigenvalues of the
Markov chain associated to $M_{W,x}$ for all $W$.

\begin{theorem}

\[ M_{G_2,x}(w) =             \left\{ \begin{array}{ll}
																																				\frac{(x+5)(x+1)}{12x^2} & \mbox{if
$d(w)=0$}\\
																																				\frac{(x+1)(x-1)}{12x^2} & \mbox{if
$d(w)=1$}\\
																																				\frac{(x-1)(x-5)}{12x^2}		& \mbox{if
$d(w)=2$}
																																				\end{array}
			\right.			 \]

\end{theorem}

\begin{proof}
	Letting $V$ be the hyperplane in $R^3$ consisting of vectors whose coordinates
add to 0, it is well known that a root system consists of $\pm (\varepsilon_i -
\varepsilon_j)$ for $i<j$ and $\pm (2 \varepsilon_i - \varepsilon_j -
\varepsilon_k)$ where $\{i,j,k\}=\{1,2,3\}$. Let $A = \varepsilon_1 -
\varepsilon_2$ and $B=-2 \varepsilon_1+\varepsilon_2+\varepsilon_3$ be a base of
positive simple roots.

	All equivalence classes $\lambda$ of subsets of $\Pi$ have size one. Some
computation gives that

\begin{eqnarray*}
e_{\emptyset} & = & \frac{1}{12} x_{\emptyset}\\
e_A & = & -\frac{1}{4} x_{\emptyset} + \frac{1}{2}x_A\\
e_B & = & -\frac{1}{4} x_{\emptyset} + \frac{1}{2}x_B\\
e_{A,B} & = & \frac{5}{12} x_{\emptyset} -\frac{1}{2}x_A - \frac{1}{2}x_B+ x_{A,B}
\end{eqnarray*}

	from which the theorem easily follows.
\end{proof}

	Let $id$ and $w_0$ be the identity and longest element of $W$. Theorems
\ref{identvalue} and \ref{longestvalue} give expressions for $M_{W,x}(id)$ and
$M_{W,x}(w_0)$. It is helpful, as
in Bergeron, Bergeron, Howlett, and Taylor
\cite{BBHT} to associate to each
$w
\in W$ an equivalence class
$\lambda$ of subsets $J$ of $\Pi$ under the action of $W$. This is done as
follows. Letting
$Fix_V(w)$ be the elements of $V$ fixed by $w$, define $A(w)=Stab_W(Fix_V(w))$.
Clearly $A(w)$ is a parabolic subgroup, conjugate to some $W_J$. Let $\lambda$ be
the equivalence class containing $J$. This $\lambda_w$ associated to $w$ will also
be called the type of $w$.

\begin{theorem} \label{identvalue} Let $m_1,\cdots,m_n$ be the exponents of
$W$. Then

\[ M_{W,x}(id) = \frac{\prod_{i=1}^n (x+m_i)}{x^n|W|} \]

\end{theorem}

\begin{proof}
	We first show that the coefficient of the identity in $e_{\lambda} = \sum_{J \in
\lambda} \frac{e_J}{|\lambda|}$ is equal to
$\frac{1}{|W|} |\{w \in W: type(w)=\lambda \}|$. Writing $e_J = \sum_{K
\subseteq J} \beta_K^J x_K$ and using the fact that the identity has coefficient 1
in each $x_K$, it is enough to show that for all $J$ of type $\lambda$,

\[ \frac{1}{|W|} |\{w \in W: type(w)= \lambda \}| = \sum_{K \subseteq J} \beta_K^J
\]

	From Bergeron, Bergeron, Howlett, and Taylor \cite{BBHT}, there is a natural map
from the descent algebra of
$W$ to the Burnside representation ring of $W$ which sends $e_J$ to $\zeta_J$ and
$x_K$ to
$Ind_{W_K}^W(1)$. Here $\zeta_J$ is the function on $W$ which takes the value $1$
if $w$ has type $\lambda$ and $0$ otherwise. This gives the equation

\[ \zeta_J = \sum_{K \subseteq J} \beta_K^J Ind_{W_K}^W(1) \]

	Take the inner product of both sides with the identity character of $W$. The
left hand side becomes $\frac{1}{|W|} |\{w \in W: type(w)= \lambda\}|$. To
evaluate the right hand side, note by Frobenius reciprocity that

\[ <Ind_{W_K}^W(1),1>_W = <1,1>_{W_K} = 1 \]

	Therefore, the right hand side becomes $\sum_{K \subseteq J} \beta_K^J$. Thus we
have shown that for all $J \in \lambda$, the coefficient of the identity in $e_J$
is equal to
$\frac{1}{|W|} |\{w \in W: type(w)=\lambda \}|$. Consequently,

\begin{eqnarray*}
M_{W,x}(id) & = & \sum_{\lambda} \frac{1}{|W|} \frac{|\{w \in W: type(w)=\lambda
\}|}{x^{\| \lambda \|}}\\
& = & \frac{1}{x^n|W|} \sum_{\lambda} x^{n-\| \lambda  \|} |\{w \in W:
type(w)=\lambda
\}|\\
& = & \frac{1}{x^n|W|} \sum_{w \in W} x^{dim(fix(w))}\\
& = & \frac{\prod_{i=1}^{n} (x+m_i)}{x^n|W|}
\end{eqnarray*}

	In the third equation, $dim(fix(w))$ is the dimension of the fixed space of $w$
in its action on $V$, the natural vector space on which $w$ acts in the relection
represenation of $W$. The third equality follows from Lemma 4.3 of Bergeron,
Bergeron, Howlett, and Taylor \cite{BBHT}, which says
that $dim(fix(w))=n-|type(w)|$. The final equality is a theorem of Shephard-Todd
\cite{ST}.
\end{proof}

\begin{theorem} \label{longestvalue} Let $m_1,\cdots,m_n$ be the exponents of
$W$. Let $w_0$ be the longest element of $W$. Then

\[ M_{W,x}(w_0) = \frac{\prod_{i=1}^n (x-m_i)}{x^n|W|} \]

\end{theorem}

\begin{proof}
	The proof is similar to that of Theorem \ref{identvalue}. It will first be shown
that the coefficient of $w_0$ in $e_{\lambda} = \sum_{J \in
\lambda} \frac{e_J}{|\lambda|}$ is equal to
$\frac{(-1)^{\| \lambda \|}}{|W|} |\{w \in W: type(w)=\lambda \}|$. Writing $e_J =
\sum_{K
\subseteq J} \beta_K^J x_K$ and using the fact that $w_0$ contributes only to
$x_{\emptyset}$, it suffices to show that for all $J$ of type $\lambda$,

\[ \beta_{\emptyset}^J = \frac{(-1)^{\| \lambda_J \|}}{|W|} |\{w \in W: type(w)=
\lambda \}| \]

	From Bergeron, Bergeron, Howlett, and Taylor \cite{BBHT}, there is a natural map
from the descent algebra of
$W$ to the Burnside representation ring of $W$ which sends $e_J$ to $\zeta_J$ and
$x_K$ to
$Ind_{W_K}^W(1)$. Here $\zeta_J$ is the function on $W$ which takes the value $1$
if $w$ has type $\lambda$ and $0$ otherwise. This gives the equation

\[ \zeta_J = \sum_{K \subseteq J} \beta_K^J Ind_{W_K}^W(1) \]

	Take the inner product of both sides with the sign character $\chi$ of $W$. The
left hand side becomes $\frac{(-1)^{\| \lambda_J \|}}{|W|} |\{w \in W: type(w)=
\lambda\}|$. To evaluate the right hand side, note by Frobenius reciprocity that

\[ <Ind_{W_K}^W(1),\chi>_W = <1,\chi>_{W_K} = \left\{ \begin{array}{ll}
																																				0 & \mbox{if
$|K|>0$}\\
		
																																				1		& \mbox{if $|K|=0$}
																																				\end{array}
			\right.		 \]

	Therefore, the right hand side becomes $\beta_{\emptyset}^J$. Thus we
have shown that for all $J \in \lambda$, the coefficient of $w_0$ in $e_J$
is equal to $\frac{(-1)^{\| \lambda_J \|}}{|W|} |\{w \in W: type(w)=\lambda \}|$.
Consequently,

\begin{eqnarray*}
M_{W,x}(w_0) & = & \sum_{\lambda} \frac{(-1)^{\| \lambda_J \|}}{|W|} \frac{|\{w
\in W: type(w)=\lambda
\}|}{x^{\| \lambda \|}}\\
& = & \frac{1}{(-x)^n|W|} \sum_{\lambda} (-x)^{n-\| \lambda  \|} |\{w \in W:
type(w)=\lambda
\}|\\
& = & \frac{1}{(-x)^n|W|} \sum_{w \in W} (-x)^{dim(fix(w))}\\
& = & \frac{\prod_{i=1}^{n} (-x+m_i)}{(-x)^n|W|}\\
& = & \frac{\prod_{i=1}^n (x-m_i)}{x^n|W|}
\end{eqnarray*}

	These equalities hold for the same reasons as in Theorem \ref{identvalue}.
\end{proof}

	Observe that left multiplication of the group algebra $Q[W]$ of $W$ by $M_{W,x}$
can be thought of as performing a random walk on $W$. The transition matrix of
this random walk is an $|W|$ by $|W|$ matrix. Theorem \ref{eigen} computes the
eigenvalues of this matrix. The eigenvalues for $W=S_n$ were determined by Hanlon
\cite{H}.

\begin{theorem} \label{eigen} The transition matrix of the random walk
arising from $M_{W,x}$ has eigenvalues $\frac{1}{x^i}$ for $0 \leq i \leq n$
with corresponding multiplicities $|{w \in W: \|type(w)\| = i}|$.
\end{theorem}

\begin{proof}
	Since the $e_{\lambda}$ decompose the identity, $Q[W] = \bigoplus_{\lambda}
e_{\lambda} Q[W]$. Since $M_{W,x} =
\sum_{\lambda}
\frac{e_{\lambda}}{x^{\| \lambda \|}}$, the eigenvalues of the action of
$M_{W,x}$ on $Q[W]$ by left multiplication are $\frac{1}{x^{\| \lambda \|}}$.
Furthermore, the eigenvalue $\frac{1}{x^i}$ occurs with mulitplicity

\[ \sum_{\lambda: \| \lambda \|=i} dim(e_{\lambda} Q[W]) \]

	Since $e_{\lambda}$ is an idempotent, $dim(e_{\lambda} Q[W])$ is the trace of
$e_{\lambda}$ regarded as a linear map from $Q[W]$ to itself. Taking the elements
of $w$ as a basis for $Q[W]$, this trace is equal to $|W|$ times the coefficient
of the identity in $e_{\lambda}$. From the proof of theorem $\ref{identvalue}$,
this coefficient of the identity in $e_{\lambda}$ is $\frac{1}{|W|} |\{w \in W:
type(w)= \lambda \}|$. This proves the theorem.
\end{proof}

{\bf Remark} In the case of the symmetric groups, Bidigare \cite{B} has computed
these eigenvalues and their mulitplicities using results of Bidigare, Hanlon, and
Rockmore \cite{BHR} on random walks arising from hyperplane arrangements. An
interesting challenge would be to prove Theorem \ref{eigen} in a similar way.

\section{Semisimple Orbits of Lie Algebras} \label{semisimple}

	This section connects the signed measures $M_{W,x}$ with semisimple
orbits of Lie algebras arising from finite groups of Lie type.

	 Let $G$ be a connected semisimple group defined over a finite field of $q$
elements. Suppose also that $G$ is simply connected. Let $\cal G$ be the Lie
algebra of $G$. Let $F$ denote both a Frobenius automorphism of $G$ and the
corresponding Frobenius automorphism of $\cal G$. Suppose that $G$ is $F$-split.
Since the derived group of $G$ is simply connected (the derived group of a simply
connected group is itself), a theorem of Springer and Steinberg \cite{SS}
implies that the centralizers of semisimple elements of $\cal G$ are connected.
Let $r$ be the rank of $G$.

	Now we define a map $\Phi$ (studied by Lehrer \cite{L} in somewhat greater
generality) from the
$F$-rational semisimple orbits $c$ of $\cal G$ to $W$, the Weyl group of $G$.
Pick $x \in {\cal G}^F \cap c$. Since the centralizers of semisimple elements of
$\cal G$ are connected, $x$ is determined up to conjugacy in $G^F$ and $C_G(x)$,
the centralizer in $G$ of $x$, is determined up to $G^F$ conjugacy. Let $T$ be a
maximally split maximal torus in $C_G(x)$. Then $T$ is an $F$-stable maximal
torus of $G$, determined up to $G^F$ conjugacy. By Proposition 3.3.3 of Carter
\cite{C}, the $G^F$ conjugacy classes of $F$-stable maximal tori of $G$ are
in bijection with conjugacy classes of $W$. Define $\Phi(c)$ to be the
corresponding conjugacy class of $W$.

	For example, in type $A_{n-1}$ the semisimple orbits $c$ of $sl(n,q)$ correspond
to monic degree $n$ polynomials $f(c)$ whose coefficient of $x^{n-1}$ vanishes.
Such a polynomial factors as $\prod_i f_i^{a_i}$ where the $f_i$ are irreducible
over $F_q$. Letting $d_i$ be the degree of $f_i$, $\Phi(c)$ is the conjugacy
class of $S_n$ corresponding to the partition $(d_i^{a_i})$.

	Two further technical concepts are helpful. As on page 28 of Carter \cite{C},
call a prime $p$ good if it divides no coefficient of any root expressed as a
linear combination of simple roots. Call a prime bad if it is not good. For
example type
$A$ has no bad primes, but
$2$ is a bad prime for type $B$. The assumption that $p$ is good will eliminate
complications involving the maximal tori of $G$ and $G^F$. Also define
$p$ to be a regular prime if the lattice of reflecting hyperplane intersections
of $W$ (including ranks of elements in the lattice) remains the same on reduction
mod
$p$. For instance, in type
$A_{n-1}$,
$p$ is not regular if $p$ divides $n$, because then $x_1=\cdots=x_n$,
$\sum {x_i}=0$ has non-trivial solutions. 	

\vspace{.5mm}

{\bf Conjecture 1:} Let $G$ be as above, and suppose that the characteristic is a
prime which is good and regular for $G$. Choose $c$ among the $q^r$
$F$-rational semisimple orbits of $\cal G$ uniformly at random. Then for all
conjugacy classes $C$ of $W$, $Prob( \Phi(c) \in C) = Prob_{M_{W,q}} (w \in C)$.
Furthermore, $M_{W,q}(w) \geq 0$ for all $w \in W$.

{\bf Remark} The assertion that $M_{W,q}(w) \geq 0$ for all $w \in W$ can be
easily checked for types $A$ and $B$ from the formulas in Section
\ref{definition} and for type $G_2$ from the formula in Section
\ref{eigenvalues}. The crucial observation (which holds for all
types), is that the bad primes for a given type are precisely those primes
which are less than the maximal exponent of $W$ but are not exponents of $W$.

	Theorems \ref{conj2id}, \ref{conj2sym}, and \ref{conj2hyp} provide evidence in
support of Conjecture 1.

\begin{theorem} \label{conj2id} Conjecture 1 holds for $G$ of all types
(i.e. $A,B,C,D,E_6, E_7,E_8,F_4,G_2$) when $C$ is the identity conjugacy class of
$W$.
\end{theorem}

\begin{proof}
	Proposition 5.9 of Lehrer \cite{L} (which uses the fact that $p$ is
regular) states that the number of
$F$-rational semisimple orbits $c$ of $\cal G$ which	satisfy $\Phi(c)=id$ is
equal to

\[ \prod_{i=1}^r \frac{q+m_i}{1+m_i} \]

	where $r$ is the rank of $G$ and $m_i$ are the exponents of $W$. Since there are
a total of $q^r$ $F$-rational semisimple orbits of $\cal G$, and because
$|W|=\prod_{i=1}^r (1+m_i)$,

	\[ Prob( \Phi(c) = id) = \frac{\prod_{i=1}^r (q+m_i)}{q^r|W|}. \]

	The proposition now follows from Theorem \ref{identvalue}.
\end{proof}

\begin{theorem} \label{conj2sym} Conjecture 1 holds for $G$ of type $A$, for all
conjugacy classes $C$ of the symmetric group $S_n$.
\end{theorem}

\begin{proof}
	Note that a monic, degree $n$ polynomial $f$ with coefficients in $F_q$ defines
a partition of $n$, and hence a conjugacy class of $S_n$, by its factorization
into irreducibles. To be precise, if $f$ factors as $\prod_i f_i^{a_i}$ where the
$f_i$ are irreducible of degree $d_i$, then $(d_i^{a_i})$ is a partition of $n$.
If the coefficient of $x^{n-1}$ in $f$ vanishes, then $f$ represents an
$F$-rational semisimple orbit $c$ of $sl(n,q)$, and the conjugacy class of
$S_n$ corresponding to the partition $(d_i^{a_i})$ is equal to $\Phi(c)$.

	Diaconis, McGrath, and Pitman \cite{DMP} have shown that if $f$ is
uniformly chosen among all monic, degree $n$ polynomials with coefficients in
$F_q$, then the measure on the conjugacy classes of $S_n$ induced by the
factorization of $f$ is equal to the measure induced by $M_{S_n,q}$. (In fact it
was this observation which led the author in the direction of Conjecture 2).

	Thus, to prove the theorem, it suffices to show that the random partition
associated to a uniformly chosen monic, degree $n$ polynomial over $F_q$ has the
same distribution as the random partition associated to a uniformly chosen monic,
degree $n$ polynomial over $F_q$ with vanishing coefficient of $x^{n-1}$. Since
the characteristic $p$ is assumed to be regular, $p$ does not divide $n$. Thus
for a suitable choice of $k$, the change of variables $x \rightarrow x+k$ gives
rise to a bijection between monic, degree $n$ polynomials with coefficient of
$x^{n-1}$ equal to $b_1$ and monic, degree $n$ polynomials with coefficient of
$x^{n-1}$ equal to $b_2$, for any $b_1$ and $b_2$. Since this bijection preserves
the partition associated to a polynomial, the theorem is proved.
\end{proof}

	Theorem \ref{conj2hyp} will confirm Conjecture 1 for all $G$ of type $B$. The
proof will use the following combinatorial objects introduced
by Reiner
\cite{R}. Let a ${\bf Z}$-word of length $m$ be a vector $(a_1,\cdots,a_m) \in
{\bf Z}^m$. For such a word define $max(a)= max(|a_i|)_{i=1}^m$. The cyclic group
$C_{2m}$ acts on
${\bf Z}$-words of length $m$ by having a generator $g$ act as $g(a_1,\cdots,a_m)
= (a_2,\cdots,a_m,-a_1)$. Call a fixed-point free orbit $P$ of this action
a primitive twisted necklace of size $m$. The group
$Z_2 \times C_m$ acts on ${\bf Z}$-words of length $m$ by having the generator $r$
of
$C_m$ act as a cyclic shift $r(a_1,\cdots,a_m)=(a_2,\cdots,a_m,a_1)$ and having
the generator $v$ of $Z_2$ act by $v(a_1,\cdots,a_m)=(-a_1,\cdots,-a_m)$. Call
a fixed-point free orbit $D$ of this action a primitive blinking necklace of size
$m$. Let a signed ornament $o$ be a set of primitive twisted necklaces and a
multiset of primitive blinking necklaces. Say that $o$ has type
$(\vec{\lambda},\vec{\mu}) =   
((\lambda_1,\lambda_2,\cdots),(\mu_1,\mu_2,\cdots))$ if it consists of
$\lambda_m$ primitive blinking neclaces of size $m$ and $\mu_m$ primitive twisted
necklaces of size $m$. Also define the size of $o$ to be the sum of the sizes of
the primitive twisted and blinking necklaces which make up $o$, and define
$max(o)$ to be the maximum of $max(D)$ and
$max(P)$ for the primitive twisted and blinking necklaces which make up $o$. 

	Reiner \cite{R} establishes the following counting lemma.

\begin{lemma} \label{Reinercount} (Reiner \cite{R})  Let $D(s,m)$ be the
number of primitive blinking necklaces $D$ such that $max(D) \leq s$. Let $P(s,m)$
be the number of primitive twisted necklaces $P$ such that $max(P) \leq s$. Then
if $q$ is an odd integer,

\[ D(\frac{q-1}{2},m) =              \left\{ \begin{array}{ll}
																																				\frac{1}{2m} \sum_{d|m \atop d \ odd} \mu(d)
(q^{\frac{m}{d}}-1) & q \geq 3, m>1\\
																																				\frac{q-1}{2} & q \geq 3, m=1\\
																																				0 & q=1 \\
																																				
																																				\end{array}
						\right. \]

\[ P(\frac{q-1}{2},m) = \left\{ \begin{array}{ll}
																																				\frac{1}{2m} \sum_{d|m \atop d \ odd} \mu(d)
	(q^{\frac{m}{d}}-1) & q \geq 3\\
																	
																																				0 & q=1 \\
																																				
																																				\end{array}
						\right. \]

\end{lemma}

	Lemma \ref{polycount} establishes an analog of Lemma \ref{Reinercount} for
special types of polynomials.

\begin{lemma} \label{polycount} Let $q$ be a positive odd integer. Let
$\tilde{I}_{m,q}$ be the number of monic, irreducible, degree
$m$ polynomials $f$ over $F_q$ satisfying $f(z)=f(-z)$. Then

\[ \tilde{I}_{2m,q} =	\left\{ \begin{array}{ll} \frac{1}{2m} \sum_{d|m \atop
d \ odd}
\mu(d)
	(q^{\frac{m}{d}}-1) & q \geq 3\\
																	
																																				0 & q=1 \\
																																				
																																				\end{array}
						\right. \]

\end{lemma}

\begin{proof}
	The case $q=1$ is clear, so assume that $q \geq 3$ is odd. Let $M_m$ be the
number of monic degree $m$ polynomials. Defining $A(t)= 1 +
\sum_{m=1}^{\infty} M_m t^m$, clearly $A(t)=\frac{1}{1-qt}$. Let $\tilde{M}_{m}$
be the number of monic degree $m$ polynomials $f$ such that $f(z)=f(-z)$.
Defining $B(t)= 1 + \sum_{m=1}^{\infty} \tilde{M}_m t^{m}$, one has that
$B(t)=\frac{1}{1-qt^2}$.  

	Observe that

	\[ A(t) = \frac{1}{1-t} \prod_{\phi: \phi(z)=\phi(-z)}
(1+\sum_{n=1}^{\infty} t^{n deg(\phi)}) \prod_{\{\phi,\tilde{\phi}\}, \phi \neq z
\atop \phi(z) \neq (-1)^{deg(\phi)} \phi(-z)} (1+\sum_{n=1}^{\infty}
t^{n deg(\phi)})^2 \]

	Here the $\phi$ are monic and irreducible, and the term $\frac{1}{1-t}$
corresponds to the contribution from the polynomial $z$. Similarly,

	\[ B(t) = \frac{1}{1-t^2} \prod_{\phi: \phi(z)=\phi(-z)}
(1+\sum_{n=1}^{\infty} t^{n deg(\phi)}) \prod_{\{\phi,\tilde{\phi}\}, \phi \neq z
\atop \phi(z) \neq (-1)^{deg(\phi)} \phi(-z)} (1+\sum_{n=1}^{\infty} t^{2n
deg(\phi)}) \]

	These equations give:

\begin{eqnarray*}
\frac{B(t)^2}{A(t^2)} & = & \frac{1}{1-t^2} \prod_{m \ even}
\frac{(1+\sum_{n=1}^{\infty} t^{mn})^{2\tilde{I}_{m,q}}}{
(1+\sum_{n=1}^{\infty} t^{2mn})^{\tilde{I}_{m,q}}}\\
& = & \frac{1}{1-t^2}
\prod_{m \ even}
\frac{(1-t^{2m})^{\tilde{I}_{m,q}}}{(1-t^{m})^{2\tilde{I}_{m,q}}}\\
& = & \frac{1}{1-t^2} \prod_{m \ even}
(\frac{1+t^{m}}{1-t^{m}})^{\tilde{I}_{m,q}}
\end{eqnarray*}

	Combining this with the explicit expressions for $A(t)$ and $B(t)$
given above shows that:

	\[ \prod_{m \ even} (\frac{1+t^m}{1-t^m})^{\tilde{I}_{m,q}} =
\frac{1-t^2}{1-qt^2} \]

	Take logarithms of both sides of this equation, using the
expansions $log(1+x)=x-\frac{x^2}{2}+\frac{x^3}{3}+\cdots$ and
$log(1-x)=-x-\frac{x^2}{2}-\frac{x^3}{3}+\cdots$.

	The left-hand side becomes:

	\[ \sum_{m \ even} \sum_{d \ odd} 2 \tilde{I}_{m,q} \frac{t^{dm}}{d} \]

	The right-hand side becomes:

	\[\sum_{m=0}^{\infty} \frac{(q^m-1) t^{2m}}{m}\]

	Comparing coefficients of $t^{2n}$ shows that $2\sum_{m|n, \frac{n}{m} \ odd}
\tilde{I}_{2m,q} m = q^n-1$. Define $L$ to be the lattice consisting of all
divisors $m$ of $n$ such that $\frac{n}{m}$ is odd. Define functions $f(m)= 2
\tilde{I}_{2m,q} m$ and $F(n)=q^n-1$ on this lattice. Moebius inversion on
this lattice implies that $f(n) = \sum_{m|n, \frac{m}{n} \ odd} \mu(m,n) F(m)$.
Thus, $\tilde{I}_{2m,q} = \frac{1}{2m} \sum_{d|m \atop d \ odd} \mu(d)
(q^{\frac{m}{d}}-1)$, as desired.
\end{proof}

	Lemma \ref{countsigned} counts the total number of signed ornaments satisfying
certain conditions. Both the result and the proof technique will be crucial in
proving Conjecture 1 for type $B$.

\begin{lemma} \label{countsigned} Let $q$ be an odd integer. The total number of
signed ornaments $o$ of size $n$ satisfying $max(o) \leq \frac{q-1}{2}$	is equal
to $q^n$.
\end{lemma}

\begin{proof}
	Let $f(z)$ be a monic polynomial over $F_q$ satisfying $f(z)=f(-z)$. Such
a polynomial can be factored uniquely as 

\[ \prod_{\{\phi_i(z),\phi_i(-z)\}} [(-1)^{deg(\phi_i)} \phi_i(z)
\phi_i(-z)]^{r_i}
\prod_{\phi_i : \phi_i(z) = \phi_i(-z)} \phi_i(z)^{s_i} \]

	where the $\phi_i$ are monic irreducible polynomials and $s_i \in \{0,1\}$.

	Hence monic polynomials satisfying $f(z)=f(-z)$ correspond to a multiset of
distinct products $(-1)^{deg(\phi)} \phi(z) \phi(-z)$ where $\phi$ is monic and
irreducible, and a set of polynomials
$\phi$ which are monic, irreducible, and satisfy $\phi(z)=\phi(-z)$. Recall that a
signed ornament corresponds to a multiset of primitive blinking necklaces and a
set of primitive twisted necklaces. Observe that there are $q^n$ monic
polynomials $f(z)$ of degree $2n$ satisfying $f(z)=f(-z)$. Lemmas
\ref{Reinercount} and \ref{polycount} show that the number of degree $2m$ monic,
irreducible polynomials satisfying $f(z)=f(-z)$ is equal to $P(\frac{q-1}{2},m)$,
the number of primitive twisted necklaces $P$ of size $m$ satisfying $max(P) \leq
\frac{q-1}{2}$.

	Thus it suffices to show that the number of distinct
products $(-1)^{deg(\phi)} \phi(z) \phi(-z)$ where $\phi$ is monic and
irreducible of degree $m$ is equal to $D(\frac{q-1}{2},m)$, the number of
primitive blinking necklaces $D$ of size $m$ satisfying $max(D) \leq
\frac{q-1}{2}$. To count the number of such products $(-1)^{deg(\phi)} \phi(z)
\phi(-z)$, note that either $\phi$ is monic, irreducible, and satisfies
$\phi(z)=(-1)^{deg(\phi)}\phi(-z)$ or else $\phi$ is monic, irreducible and does
not satisfy $\phi(z)=(-1)^{deg(\phi)} \phi(-z)$, this latter case arising in two
possible ways. Thus the number of such products is equal to
$\frac{A(m,q)+B(m,q)}{2}$, where $A(m,q)$ is the number of monic, irreducible
$\phi$ of degree $m$ satisfying $\phi(z)=\phi(-z)$, and $B(m,q)$ is the the number
of monic, irreducible $\phi$ of degree $m$. Lemma \ref{polycount} shows that

\[ A(m,q) =	\left\{ \begin{array}{ll} \frac{1}{m} \sum_{d|m \atop
d \ odd}
\mu(d)
	(q^{\frac{m}{2d}}-1) & m \ even\\
																	
																																				0 & m \ odd\\
																																				
																																				\end{array}
						\right. \]

It is well known that $B(m,q) = \frac{1}{m} \sum_{d|m} \mu(d) q^{\frac{m}{d}}$.
Easy manipulations show that

\[ \frac{A(m,q)+B(m,q)}{2} =              \left\{ \begin{array}{ll}
																																				\frac{1}{2m} \sum_{d|m \atop d \ odd} \mu(d)
(q^{\frac{m}{d}}-1) & q \geq 3, m>1\\
																																				\frac{q-1}{2} & q \geq 3, m=1\\
																																				0 & q=1 \\
																																				
																																				\end{array}
						\right. \]

	Thus $\frac{A(m,q)+B(m,q)}{2} = D(\frac{q-1}{2},m)$, and the lemma is proved.
\end{proof}

	With these lemmas in hand, Conjecture 1 can be proved for type $B$.

\begin{theorem} \label{conj2hyp} Conjecture 1 holds for $G$ of type $B$, for all
conjugacy classes $C$ of the hyperoctahedral group $B_n$.
\end{theorem}

\begin{proof}
	Note that because $2$ is a bad prime for type $B$, it can be assumed that the
characteristic is odd.

	Recall that the type of a signed ornament is parameterized by pairs of vectors
$(\vec{\lambda},\vec{\mu})$, where $\lambda_i$ is the number of primitive
blinking necklaces of size $i$ and $\mu_i$ is the number of primitive twisted
necklaces of size $i$. It is well known from the theory of wreath products that
the conjugacy classes of the hyperoctahedral group $B_n$ are also parameterized
by pairs of vectors $(\vec{\lambda},\vec{\mu})$, where $\lambda_i(w)$ and
$\mu_i(w)$ are the number of positive and negative cycles of $w \in B_n$
respectively.

	The first step of the proof will be to show that the measure induced on pairs
$(\vec{\lambda},\vec{\mu})$ by choosing a random signed ornament $o$ of size $n$
satisfying $max(o) \leq \frac{q-1}{2}$ is equal to the measure induced on
pairs $(\vec{\lambda},\vec{\mu})$ by choosing $w \in B_n$ according to the measure
$M_{B_n,q}$ and then looking at its conjugacy class.

	From the definition of descents given in Section \ref{definition}, it is easy to
see that if one introduces the following linear order $\Lambda$ on the set of
non-zero integers:

	\[ +1 <_{\Lambda} +2 <_{\Lambda} \cdots +n<_{\Lambda} \cdots <_{\Lambda} -n
<_{\Lambda} \cdots <_{\Lambda} -2 <_{\Lambda} -1 \]

then $d(w)$, the number of descents of $w \in B_n$, can be defined as $|\{i: 1
\leq i \leq n: w(i) <_{\Lambda} w(i+1)\}|$. Here $w(n+1)=n+1$ by convention.

	Reiner \cite{R} proves that there is a bijection between signed ornaments
$o$ of size $n$ satisfying $max(o) \leq \frac{q-1}{2}$ and pairs $(w,\vec{s})$
where
$w
\in B_n$ and $\vec{s}=(s_1,\cdots,s_n) \in {\bf N}^n$ satisifies $\frac{q-1}{2}
\geq s_1 \geq \cdots \geq s_n \geq 0$ and $s_i > s_{i+1}$ when $w(i) <_{\Lambda}
w(i+1)$ (i.e. when $w$ has a descent at position $i$). Further, he shows that
the type of $o$ is equal to the conjugacy class vector of $w$.

	It is easy to see that if $w$ has $d(w)$ descents, then the number of
$\vec{s}$ such that $\frac{q-1}{2}
\geq s_1 \geq \cdots \geq s_n \geq 0$ and $s_i > s_{i+1}$ when $w(i) <_{\Lambda}
w(i+1)$ is equal to

\[ {\frac{q-1}{2}+n-d(w) \choose n} = \frac{(q+1-2d(\pi)) \cdots
(q+2n-1-2d(\pi))}{2^n n!} \]

	Lemma \ref{polycount} shows that there are $q^n$ signed ornaments $f$ of size
$n$ satisfying $max(f) \leq \frac{q-1}{2}$. Thus we
conclude that choosing a random signed ornament induces a measure on $w \in
B_n$ with mass on $w$ equal to

\[ \frac{(q+1-2d(\pi)) \cdots (q+2n-1-2d(\pi))}{q^n |B_n|} \]

	By the remarks in Section \ref{definition}, this is exactly the mass on $w$
under the measure $M_{B_n,q}$. Since in Reiner's bijection the type of $o$ is
equal to the conjugacy class vector of $w$, we have proved that the measure on
conjugacy classes $(\vec{\lambda},\vec{\mu})$ of $B_n$ induced by choosing $w$
according to $M_{B_n,q}$ is equal to the measure on conjugacy
classes $(\vec{\lambda},\vec{\mu})$ of $B_n$ induced by choosing a signed
ornament uniformly at random and taking its type.

	The second step in the proof is to show that if $f$ is chosen uniformly
among the $q^n$ semisimple orbits of $o(2l+1,q)$, then the chance that $\Phi(f)$
is the conjugacy class $(\vec{\lambda},\vec{\mu})$ of $B_n$ is equal to the
chance that a signed ornament chosen randomly among the $q^n$ signed ornaments $o$
of size $n$ satisfying $max(o) \leq \frac{q-1}{2}$ has type
$(\vec{\lambda},\vec{\mu})$.

	It is well known that the semisimple orbits of $Spin(2n+1,q)$ on $o(2l+1,q)$
correspond to monic, degree $2n$ polynomials $f$ satisfying $f(z)=f(-z)$. It is
also not difficult to see that $\Phi(f)$ can be described as follows. Factor $f$
uniquely as 

\[ \prod_{\{\phi_i(z), \phi_i(-z)\}} [(-1)^{deg(\phi_i)} \phi_i(z)
\phi_i(-z)]^{r_i}
\prod_{\phi_i : \phi_i(z) = \phi_i(-z)} \phi_i(z)^{s_i} \]

	where the $\phi_i$ are monic irreducible polynomials and $s_i \in \{0,1\}$. Then
let $\lambda_i(f)=\sum r_i$ and $\mu_i(f)=\sum s_i$. Lemmas
\ref{Reinercount} and \ref{polycount} show that the number of degree $2m$ monic,
irreducible polynomials satisfying $f(z)=f(-z)$ is equal to
the number of primitive twisted necklaces $P$ of size $m$ satisfying $max(P) \leq
\frac{q-1}{2}$. Lemma \ref{polycount} shows that the number of distinct
products $(-1)^{deg(\phi)} \phi(z) \phi(-z)$ where $\phi$ is monic and
irreducible of degree $m$ is equal to the number of primitive blinking necklaces
$D$ of size $m$ satisfying $max(D) \leq \frac{q-1}{2}$. This proves the theorem.
\end{proof}

{\bf Remarks}

\begin{enumerate}

\item It is worth pointing out that Conjecture 1 would be false if instead of
choosing $c$ uniformly among the $q^r$ $F$-rational semisimple orbits of ${\cal
G}$, $c$ were chosen uniformly among the $q^r$ semisimple conjugacy classes of
$G^F$. For a simple counterexample, take $G=SL_3(5)$ and $C$ the identity
conjugacy class of $S_3$. There are only five monic polynomials $f$ with
coefficients in $F_5$ which factor into linear terms and satisfy $f(0)=1$. The
analog of Conjecture 2 would predict that there are seven.		

\item Let $Q_{p'}$ be the additive group of rational numbers of the form
$r/s$ where $r,s$ are integers and $s$ is not divisible by $p$. Conjecture 1 leads
us to speculate that after fixing some extra structure such as a Borel subgroup
and an isomorphism between
$\bar{F}_q^*$ and
$Q_{p'}/Z$, there should be a canonical way to associate to an
$F$-rational semisimple orbit $c$ of ${\cal G}$ an {\it element} $w$ of $W$,
inducing the measure $M_{W,q}$ on $W$. Furthermore, the conjugacy class of $w$
should be equal to $\Phi(c)$.

\end{enumerate}

\section{Acknowledgements} The author thanks Persi Diaconis for many useful
references. Dick Gross has also been helpful in clarifying some points about
algebraic groups.

\end{document}